# General Elementary Direct Proof of Fermat's Last Theorem

Hua Jiang

*Department of Chemical & Environmental Engineering, University of Arizona, Tucson, AZ*
Email: huaj@email.arizona.edu

**Abstract**

This paper presents a novel direct elementary proof for Fermat's Last Theorem. We use algebra, modular math, and binomial series to develop inherent mathematical relationships hidden within Fermat's Last Theorem. With these derived relationships, we are able to develop general pattern applicable for all positive integers of $n$. Finally, we are able to confirm and complete the direct proof for Fermat's Diophantine equation for all $n$.

**Introduction**

Around 1630, the Fermat's Last Theorem (FLT) was noted by Pierre de Fermat in his copy of Diophantus's Arithmetica [1]. Fermat famously wrote that he had a proof that was too large to fit in the margin. It was not until 1995 that a successful proof by contradiction was published by Andrew Wiles using many modern mathematical techniques developed since Fermat [2]. However, only a few mathematicians are knowledgeable and able to comprehend this proof. Questions have risen whether there exists a simpler proof for FLT, especially since mathematical knowledge and understanding from Fermat time has expanded. Most attempted FLT proofs and Wiles' proof have relied on the proof by contradiction. In this proof, we present a general elementary direct proof of FLT for all positive $n$ integers.

To prove Fermat's last theorem (FLT), we first utilize this statement defined by FLT [1]:

$$x^n + y^n = z^n.$$

where $x$, $y$, and $z$ and $n$ are positive integers. We will establish inherent relationships and develop a general direct proof using basic algebra and modular math to prove the claim that there no integer solutions of FLT for $n > 2$ and there are integer solutions for $n = 1$ and 2.

**Definition 1.** Fermat's Last Theorem is redefine as

$$x^n + (x + \beta)^n = (x + \beta + \alpha)^n \tag{1}$$

where $\beta = y - x$, $\alpha = z - y$, $\beta \neq 0$, and $\alpha > 0$.

**Lemma 1.** Let (1) be represented by two rational numbers, $k$ and $l$, and one integer, $n$:

$$k^n + l^n = (l + 1)^n \tag{2}$$

where

$$l = k + \frac{\beta}{\alpha}. \tag{3}$$

*Proof.* By applying a modulus of $\alpha$ on equation (1), $x^n$ is found to be divisible by $\alpha$:

$$x^n = 0 \bmod \alpha, \tag{4}$$

$$x \propto \alpha.$$

We can claim that $x$ is proportional to $\alpha$. Let $k$, a rational number, be the proportional relationship between $x$ and $\alpha$:

$$x = k\alpha. \tag{5}$$

Variable $k$ cannot represent all rational numbers. Values of $k$ that does not satisfy (4) are invalid to FLT.

For example, when $x = 42, \alpha = 24$, (4) is satisfied for $n > 2$:

$$x^1 = 18 \bmod 24, x^2 = 12 \bmod 24, x^{n>2} = 0 \bmod 24.$$

For example, when $x = 42, \alpha = 36$, (4) is satisfied for $n > 1$:

$$x^1 = 6 \bmod 36, x^{n>1} = 0 \bmod 36.$$

Applying equation (5) into equation (1) results in

$$(k\alpha)^n + (k\alpha + \beta)^n = (k\alpha + \beta + \alpha)^n.$$

Divide by $\alpha^n$:

$$k^n + \left(k + \frac{\beta}{\alpha}\right)^n = \left(k + \frac{\beta}{\alpha} + 1\right)^n.$$

Apply substitution from (3) to get (2).

**Lemma 2.** Equation (2) can be restructured as

$$fy = \alpha^n ((v+1)^n - 1). \tag{6}$$

where $v$ is a rational number, $n$ as an integer, $y$ as an integer, $f$ as modulus multiple, and

$$k = v + 1. \tag{7}$$

*Proof*: Using the binomial theorem, expand equation (2) into a sum series:

$$k^n = (v+1)^n = (l+1)^n - l^n,$$

$$\sum_{i=0}^{n} \frac{n!}{i!(n-i)!} v^i = \sum_{i=0}^{n-1} \frac{n!}{i!(n-i)!} l^i.$$

At $i = 0$, the $0^{th}$ term for both sum series is 1. This cancels out and leaves factors $v$ and $l$ respectively for the left and right side. Let $t$, a rational number, be the proportional relationship between $v$ and $l$:

$$l = vt. \tag{8}$$

Substitute in $l$ and merge the sum series:

$$\sum_{i=0}^{n} \frac{n!}{i!\,(n-i)!} v^i = \sum_{i=0}^{n-1} \frac{n!}{i!\,(n-i)!} v^i t^i.$$

Take the $n^{\text{th}}$ term out and merge the sum series:

$$v^n + \sum_{i=1}^{n-1} \frac{n!}{i!\,(n-i)!} v^i = \sum_{i=1}^{n-1} \frac{n!}{i!\,(n-i)!} v^i t^i,$$

$$v^n = \sum_{i=1}^{n-1} \frac{n!}{i!\,(n-i)!} (t^i - 1) v^i.$$

There is another common factor, $t - 1$, in the sum series. Let $g$, a rational number, be the proportional relationship between $v$ and $t - 1$.

$$v = g(t - 1) \tag{9}$$

Define $l$ in terms of $v$ and $g$ with equations (8) and (9).

$$l = vt = \frac{v^2}{g} + v$$

Calculate the boundary for $g$:

$$l = \frac{v^2}{g} + v = k + \frac{\beta}{\alpha},\, k = v + 1,$$

$$\frac{v^2}{g} = 1 + \frac{\beta}{\alpha},$$

$$\frac{v^2}{g} = \frac{\alpha + \beta}{\alpha} = \zeta, \tag{10}$$

where $\zeta$ is a positive rational constant for a set of solutions. For example, the multiples of Pythagoras triple set (3, 4, 5) has a $\zeta$ of 2 and the multiples of Pythagoras triple set (5, 12, 13) has a $\zeta$ of 8.

From (10), $g$ can be defined as a positive rational number.

$$g = \frac{\alpha}{\alpha + \beta} v^2,$$

$$\alpha + \beta > 0,\, \alpha > 0,$$

$$0 < g < \infty.$$

Define $y$ in terms of $v$, $g$, and $\alpha$ with equations (3), (5), and $l$:

$$\frac{v^2}{g} + v = k + \frac{\beta}{\alpha},$$

$$y = k\alpha + \beta = \left(\frac{v^2}{g} + v\right)\alpha. \tag{11}$$

Substitute (5) in (1)

$$(k\alpha)^n = (k\alpha + \beta + \alpha)^n - (k\alpha + \beta)^n.$$

Apply substitutions with (7) and $k\alpha + \beta$ to the above equation and factor out $\alpha^n$.

$$\alpha^n(v+1)^n = \alpha^n\left(\frac{v^2}{g} + v + 1\right)^n - \alpha^n\left(\frac{v^2}{g} + v\right)^n \tag{12}$$

Apply a modulus of integer $\left(\frac{v^2}{g} + v\right)\alpha$ to (12):

$$\alpha^n(v+1)^n = \alpha^n \bmod \left(\frac{v^2}{g} + v\right)\alpha,$$

$$\alpha^n((v+1)^n - 1) = f\left(\frac{v^2}{g} + v\right)\alpha = fy.$$

**Lemma 3.** From equation (6), $fy$ can be derived in an alternative binomial form.

$$fy = \sum_{i=1}^{n-1} \binom{n}{i} \alpha^i y^{n-i} \tag{13}$$

*Proof*: Equation (6) can be simplified with (7) and then by (5):

$$fy = \alpha^n(k^n - 1),$$

$$fy = x^n - \alpha^n.$$

By substituting $x^n$ defined by (1), we find another representation of $fy$:

$$x^n + y^n = (y + \alpha)^n,$$

$$fy = (y + \alpha)^n - y^n - \alpha^n.$$

**Proof of FLT**

We have now generated inherent hidden relationships within the FLT equation from the definitions of $x, y, z$ and $n$ as described in FLT. To confirm the validity of FLT for all $n$, we first equate the two $fy$ equations (6) and (13) and obtain a polynomial of $v$:

$$\alpha^n((v+1)^n - 1) = \sum_{i=1}^{n-1} \binom{n}{i} \alpha^i y^{n-i}.$$

Obtain $y$ in terms of $\alpha$, $\zeta$, and $v$ with (10) and (11)

$$y = (\zeta + v)\alpha.$$

Substitute in $y$ and cancel $\alpha^n$ on both sides

$$(v+1)^n - 1 = \sum_{i=1}^{n-1} \binom{n}{i} (\zeta + v)^{n-i}. \tag{14}$$

To develop a general pattern, we will investigate values of $n$ from 1 to 5.

**Case $n = 1$**: The FLT equation becomes very simple and does not require the application of the binomial expansion, so we use (2) when $n = 1$:

$$k^1 = (l+1)^1 - l^1 = 1.$$

We find that $k$ is equal to 1 and $l$ is independent of $k$. This results in infinite solutions for $n = 1$.

**Case $n = 2$**: Set $n$ equal to 2 in (14) and expand

$$(v+1)^2 - 1 = \sum_{i=1}^{1} \binom{2}{i}(\zeta+v)^{2-i} = 2(\zeta + v).$$

Apply the substitution for $\zeta$ with (10) and expand

$$v^2 + 2v = g\zeta + 2v = 2\zeta + 2v.$$

By equating coefficients, the equation reduces to:

$v^1$: $\qquad\qquad\qquad\qquad 2v = 2v,$

$v^0$: $\qquad\qquad\qquad\qquad g\zeta = 2\zeta.$

We find that $g$ is equal to 2 irrespective of $v$. For both equations, $v$ is not limited to any additional condition. Thus, there are infinite solutions for $n = 2$.

**Case $n = 3$**: Set $n$ equal to 3 in (14) and expand:

$$(v+1)^3 - 1 = \sum_{i=1}^{2} \binom{3}{i}(\zeta+v)^{3-i},$$

$$v^3 + 3v^2 + 3v = \binom{3}{1}(\zeta+v)^2 + \binom{3}{2}(\zeta+v)^1,$$

$$v^3 + 3v^2 + 3v = \binom{3}{1}(\zeta^2 + 2v\zeta + v^2) + \binom{3}{2}(\zeta + v).$$

After common terms are cancelled out, the equation becomes

$$v^3 = \binom{3}{1}(\zeta^2 + 2v\zeta) + \binom{3}{2}\zeta.$$

We apply substitutions for factors of $v^2$ with (10) to form a first order polynomial of $v$:

$$g\zeta v = 3\zeta^2 + 6\zeta v + 3\zeta.$$

Equate the coefficients for $v^1$ and $v^0$:

$v^1$: $\qquad\qquad\qquad\qquad g\zeta v = 6\zeta v,$

$v^0$:
$$3\zeta^2 + 3\zeta = \sum_{i=1}^{2} \binom{3}{i}\zeta^i = 0.$$

We obtain a value of 6 for $g$ and $-1$ for $\zeta$. However, $\zeta$ cannot be negative. Thus, there are no rational solutions for $n = 3$.

**Case $n = 4$:** Set $n$ equal to 4 in (14) and expand:

$$(v+1)^4 - 1 = \sum_{i=1}^{3} \binom{4}{i}(\zeta+v)^{4-i},$$

$$v^4 + 4v^3 + 6v^2 + 4v = \binom{4}{1}(\zeta+v)^3 + \binom{4}{2}(\zeta+v)^2 + \binom{4}{3}(\zeta+v)^1,$$

$$v^4 = \binom{4}{1}(\zeta^3 + 3\zeta^2 v + 3\zeta v^2) + \binom{4}{2}(\zeta^2 + 2v\zeta) + \binom{4}{3}(\zeta)^1.$$

Two binomial series can be formed from the first term from the three groups and the second term from the first two groups:

$$4\zeta^3 + 6\zeta^2 + 4\zeta = \sum_{i=1}^{3}\binom{4}{i}\zeta^i,$$

$$12\zeta^2 v + 12\zeta v = 4v\sum_{i=1}^{2}\binom{3}{i}\zeta^i.$$

After common terms are cancelled out, the equation becomes

$$v^4 = 12\zeta v^2 + \sum_{i=1}^{3}\binom{4}{i}\zeta^i + 4v\sum_{i=1}^{2}\binom{3}{i}\zeta^i.$$

We apply substitutions for factors of $v^2$ with (10) to form a first order polynomial of $v$:

$$(g\zeta)^2 = 12g\zeta^2 + \sum_{i=1}^{3}\binom{4}{i}\zeta^i + 4v\sum_{i=1}^{2}\binom{3}{i}\zeta^i.$$

Equate the coefficients for $v^1$ and $v^0$:

$v^1$:
$$4v\sum_{i=1}^{2}\binom{3}{i}\zeta^i = 12\zeta^2 v + 12\zeta v = 0,$$

$v^0$:
$$g^2\zeta^2 = 12g\zeta^2 + \sum_{i=1}^{3}\binom{4}{i}\zeta^i$$

From equating the coefficients for $v^1$, we obtain a value of $-1$ for $\zeta$. However, $\zeta$ cannot be negative. Thus, there are no rational solutions for $n = 4$.

**Case $n = 5$:** Set $n$ equal to 5 in (14) and expand:

$$(v+1)^5 - 1 = \sum_{i=1}^{4} \binom{5}{i}(\zeta + v)^{5-i},$$

$$v^5 + 5v^4 + 10v^3 + 10v^2 + 5v = \binom{5}{1}(\zeta+v)^4 + \binom{5}{2}(\zeta+v)^3 + \binom{5}{3}(\zeta+v)^2 + \binom{5}{4}(\zeta+v)^1,$$

$$v^5 = \binom{5}{1}(\zeta^4 + 4\zeta^3 v + 6\zeta^2 v^2 + 4\zeta v^3) + \binom{5}{2}(\zeta^3 + 3\zeta^2 v + 3\zeta v^2) + \binom{5}{3}(\zeta^2 + 2v\zeta)$$
$$+ \binom{5}{4}(\zeta)^1.$$

After common terms are cancelled out and collected terms into a binomial series, the equation becomes:

$$v^5 = 30\zeta v^3 + \sum_{i=1}^{5-1} \binom{5}{i}\zeta^i + 5v\sum_{i=1}^{4-1}\binom{4}{i}\zeta^i + 10v^2\sum_{i=1}^{3-1}\binom{3}{i}\zeta^i.$$

We see there is a developing pattern with the binomials:

$$v^5 = 30\zeta v^3 + \binom{5}{0}\sum_{i=1}^{4}\binom{5}{i}\zeta^i + \binom{5}{1}v\sum_{i=1}^{3}\binom{4}{i}\zeta^i + \binom{5}{2}v^2\sum_{i=1}^{2}\binom{3}{i}\zeta^i.$$

Apply (10) for factors of $v^2$ and equate the coefficients for $v^1$ and $v^0$:

$v^1$:
$$(g\zeta)^2 v = 30\zeta(g\zeta)v + \binom{5}{1}v\sum_{i=1}^{3}\binom{4}{i}\zeta^i,$$

$v^0$:
$$\binom{5}{0}\sum_{i=1}^{4}\binom{5}{i}\zeta^i + \binom{5}{2}(g\zeta)\sum_{i=1}^{2}\binom{3}{i}\zeta^i = 0.$$

From equating coefficients for $v^0$, we get a polynomial of $g$ and $\zeta$ equal to zero. However, since the coefficients of the polynomial are positive and, $g$ and $\zeta$ are limited to positive rational numbers, then the equation cannot be satisfied. Thus, there are no solutions for $n = 5$.

We developed a generalized equation from equating (6) and (13) from the previous cases:

$$v^n = n(n-1)\zeta v^{n-2} + \binom{n}{0}\sum_{i=1}^{n-1}\binom{n}{i}\zeta^i + \binom{n}{1}v\sum_{i=1}^{n-2}\binom{n-1}{i}\zeta^i + \binom{n}{2}v^2\sum_{i=1}^{n-3}\binom{n-2}{i}\zeta^i + \cdots$$
$$+ \binom{n}{n-3}v^{n-3}\sum_{i=1}^{2}\binom{3}{i}\zeta^i.$$

As shown in case $n = 3$, 4, and 5, we can group powers of even and odd powers of $v$ into two equations. This is due to the substitution with (10) that reduces even and odd powers of $v$ respectively to $v^0$ and $v^1$. From equating the coefficients, there will be two equations represented as a sum of a group of binomial series which one equation is equal to zero and the other is equal to a function of $g$ and $\zeta$.

For odd $n > 2$, the zero equation is obtained when we equate the coefficients for $v^0$:

$$\binom{n}{0} \sum_{i=1}^{n-1} \binom{n}{i} \zeta^i + \binom{n}{2} v^2 \sum_{i=1}^{n-3} \binom{n-2}{i} \zeta^i + \cdots + \binom{n}{n-3} v^{n-3} \sum_{i=1}^{2} \binom{3}{i} \zeta^i = 0.$$

Apply (10) for $v^2$ and formulate a double sum series

$$\sum_{j=0}^{\frac{n-1}{2}-1} \binom{n}{2j} (g\zeta)^j \sum_{i=1}^{n-2j-1} \binom{n-2j}{i} \zeta^i = 0. \tag{15}$$

For even $n > 2$, the zero equation is obtained when we equate the coefficients for $v^1$:

$$\binom{n}{1} v \sum_{i=1}^{n-2} \binom{n-1}{i} \zeta^i + \binom{n}{3} v^3 \sum_{i=1}^{n-4} \binom{n-3}{i} \zeta^i + \cdots + \binom{n}{n-4} v^{n-4} \sum_{i=1}^{2} \binom{3}{i} \zeta^i = 0.$$

Apply (10) for $v^2$ and formulate a double sum series

$$v \sum_{j=0}^{\frac{n}{2}-2} \binom{n}{2j+1} (g\zeta)^j \sum_{i=1}^{n-2j-2} \binom{n-2j-1}{i} \zeta^i = 0. \tag{16}$$

For both zero equations (15) and (16), the coefficients of the polynomial are all positive. Since $g$ and $\zeta$ are positive rational numbers, the equations can never be satisfied to equal zero. Thus there are no rational solutions for $n > 2$. With the proven cases of $n = 1$ and $n = 2$, this completes the general elementary direct proof of FLT for positive integers of $n$.

**Acknowledgements**

The author wish to thank Society for Industrial and Applied Mathematics (SIAM) Student Chapter at the University of Arizona for their critique and discussion. The author also wish to thank Dr. Vasiliki Karanikola, Dr. Eduardo Saez, and the University of Arizona for the support.